\documentclass[11pt,a4paper]{article}
\usepackage[english]{babel}

\usepackage{graphicx}
\begin{document}

\title{Identification of Convection Heat Transfer Coefficient \\ of Secondary Cooling Zone of CCM \\ based on Least Squares Method and Stochastic Approximation Method}

\author{G. O. Ivanova \\ IAMM NAN of Ukraine, Donetsk\\ivanova@iamm.ac.donetsk.ua}

\maketitle

\begin{abstract}

The detailed mathematical model of heat and mass transfer of steel
ingot of curvilinear continuous casting machine is proposed.
 The process of heat and mass transfer is described by nonlinear partial
differential equations of parabolic type. Position of phase boundary
is determined by Stefan conditions. The temperature of cooling water
in mould channel is described by a special balance equation.
Boundary conditions of secondary cooling zone include radiant and
convective components of heat exchange and account for the complex
mechanism of heat-conducting due to airmist cooling using compressed
air and water. Convective heat-transfer coefficient of secondary
cooling zone is unknown and considered as distributed parameter. To
solve this problem the algorithm of initial adjustment of parameter
and the algorithm of operative adjustment are developed.

\end{abstract}

\section{Introduction}

Improved computing significantly increased
 role of mathematical
modeling in research of thermo-physical processes. This, in turn,
imposes stricter requirements towards accuracy and efficiency of
mathematical models.

It is well known that successful modeling mostly depends on the
right choice of a model, which is directly affected by reliability
of thermo-physical parameters used. Frequently, empirical data alone
can not provide sufficient information about one-valuedness
conditions.

Therefore recently the big attention is given to the solution of
inverse problems of heat conduction, in which it is necessary to
define thermophysical properties of an object on available
(frequently rather limited) information about temperature field. In
particular thus it is possible to identify boundary conditions.
There are difficulties in choice of some parameters of process for
development of mathematical models of technological processes.

While modeling process for specific industrial conditions it is
necessary to determine some thermal or physical parameters each
time, in particular convective heat-transfer coefficient (CHTC) on a
surface of an ingot in the secondary cooling zone which depends on
many factors. It is connected by that the convective heat transfer
coefficient value is influenced with set of various factors.
Besides, CHTC value can vary strongly enough in a time and on space
coordinates. Thus, there is a problem of identification of the CHTC
as distributed parameter.

In the given work algorithms of initial adjustment of parameter when
at the disposal of there is enough plenty of points in which the
temperature on a surface of an ingot is measured, and operative
adjustment when the temperature is measured only in one point on a
surface are considered.

\section{Statement of problem} The thermal field of the moving steel
ingot and mold wall in the system of coordinates attached to
motionless construction of CCM is considered~[1]. In fig. 1 the
diagram of CCM is introduced.

\begin{figure}[ht]
\centering

\includegraphics[width=140mm, clip]{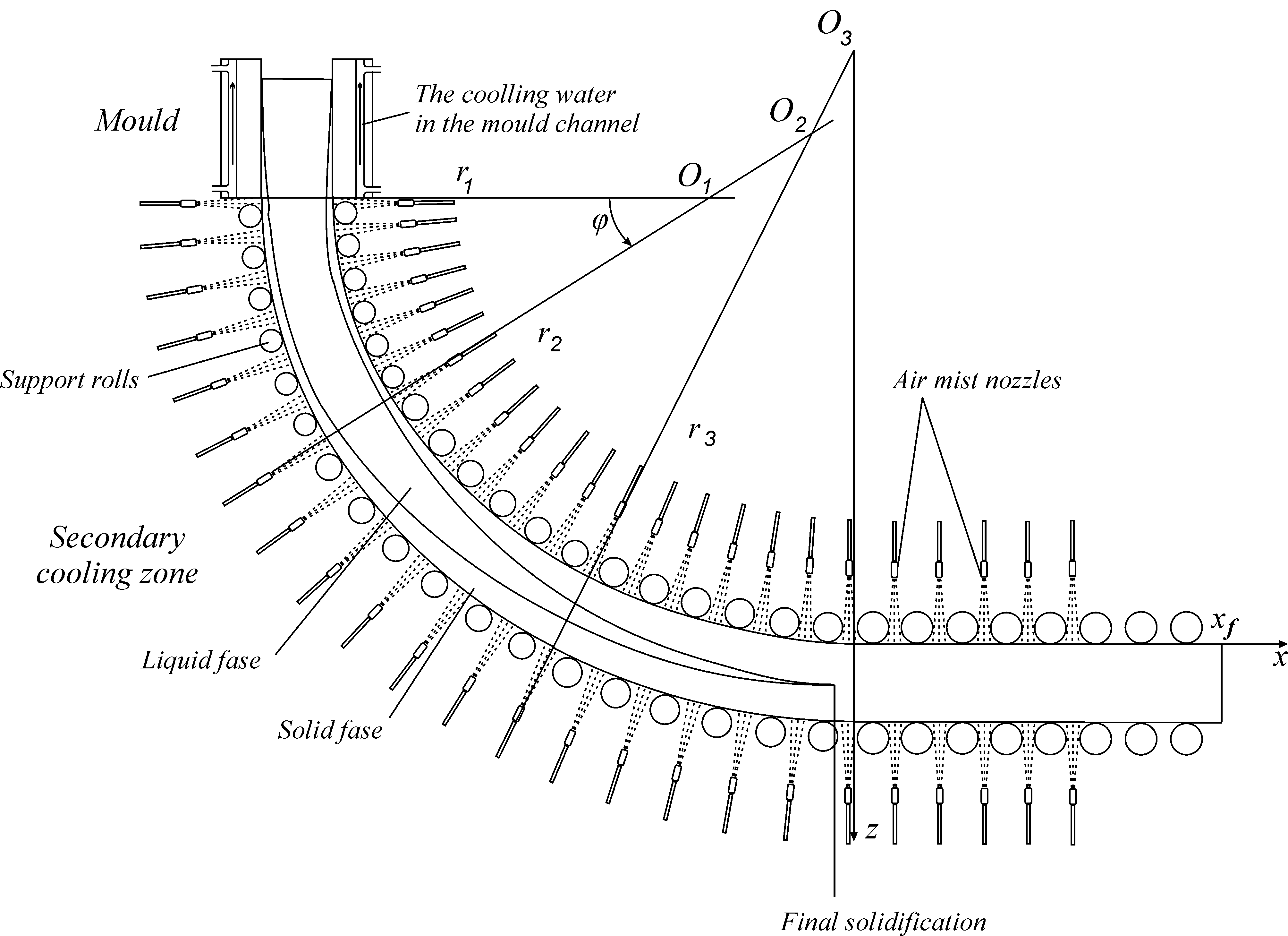}
\caption{}
\end{figure}

The heat conduction in the steel ingot in the mold area is described
by nonstationary, nonlinear heat and mass transfer equation:

$$
\begin{array}{l}
 {\displaystyle \frac{\partial T(\tau ,x,z)}{\partial \tau } + v(\tau
)\frac{\partial T\left( {\tau ,x,z} \right)}{\partial z} = } \\ \\
 {\displaystyle=\frac{1}{c(T,x,z)\rho (T,x,z)}\left\{ {\frac{\partial }{\partial
x}\left[ {\lambda (T,x,z)\frac{\partial T}{\partial x}} \right] +
\frac{\partial }{\partial z}\left[ {\lambda (T,x,z)\frac{\partial
T}{\partial z}} \right]} \right\},}
\\ \\
{\displaystyle\ \qquad \qquad \qquad \qquad \qquad 0 < x < l, \quad
0 < z < Z}
\end{array} \eqno{(1)}
$$

\noindent and the boundary conditions:

$$
\begin{array}{c}
{\displaystyle\ - \lambda (T,x)\frac{\partial T}{\partial z} = 0,
\quad 0 \le x \le l, }
\\ \\
{\displaystyle\  \left. {\frac{\partial T}{\partial x}} \right|_{x =
0} = 0, \quad 0 \le z \le Z, }
\\ \\
{\displaystyle\ \lambda (T,z)\left. {\frac{\partial T}{\partial x}}
\right|_{x = l} = \frac{\lambda _{gz} }{\delta }\left( {\left.
T\right|_{x = l + \delta } - \left. T\right|_{x = l} } \right) +
\sigma_n\left[ {\left( {\frac{\left. T\right|_{x = l + \delta }
}{100}} \right)^4 - \left( {\frac{\left. T\right|_{x = l} }{100}}
\right)^4} \right],}
\\ \\
\quad 0 \le z \le Z,
\end{array} \eqno{(2)}
$$

\noindent where $v(\tau )$ -- withdrawal rate, $2l$ -- ingot
thickness, $Z$ -- height of ingot in the mould, $T(\tau ,x,z)$ --
metal temperature, $c(T,x,z)$ -- metal specific heat, $\rho (T,x,z)$
-- density, $\lambda (T,x,z)$ -- thermal conduction, $\delta ~$ --
effective thickness of air gap between ingot and the mould wall,
$\lambda _{gz} $ -- thermal conduction coefficient of gap gas
mixture, $\left. T \right|_{x = l} $ -- surface temperature of the
ingot, $\left. T \right|_{x = l + \delta } $ -- surface temperature
of mold wall, $\sigma_n$ -- the resulted radiation coefficient.

Conditions of equality of temperatures and Stefan conditions, and
also boundary and initial conditions for the phase boundary are set:

$$
\begin{array}{c}
{\displaystyle\ T\left. {(\tau ,x,z)} \right|_{x = \xi _ - (\tau
,z)} = \left. {T(\tau ,x,z)} \right|_{x = \xi _ + (\tau ,z)} =
T_{kr} ,}
\\ \\
{\displaystyle\ \left. {\lambda (T,x,z)\frac{\partial T}{\partial
\bar {n}}} \right|_{x = \xi _ - (\tau ,z)} - \left. {\lambda
(T,x,z)\frac{\partial T}{\partial \bar {n}}} \right|_{x = \xi _ +
(\tau ,z)} = \mu \rho (T_{kr} )\left( {\frac{\partial \xi }{\partial
\tau } + v \cdot \frac{\partial \xi }{\partial z}} \right), }
\\ \\
{\displaystyle\ \quad 0 \le z \le Z, }
\\ \\
{\displaystyle\ \xi (\tau ,0) = l, \quad \xi (0,z) = \xi _0 (z),}
\end{array} \eqno{(3)}
$$

\noindent where $\xi $ -- the phase boundary function of two
variables $x = \xi (\tau ,z)$, $\mu $ -- crystallization latent
heat, $T_{kr} $ -- crystallization temperature (average of the
interval ``liquidus -- solidus''), $\bar {n}$ -- normal to the
boundary of phases.

Heat equation for mould walls:

$$
\begin{array}{c}
{\displaystyle\ \frac{\partial T(\tau ,x,z)}{\partial \tau } =
\frac{1}{c(T,x,z)\rho (T,x,z)}\left\{ {\frac{\partial }{\partial
x}\left[ {\lambda (T,x,z)\frac{\partial T}{\partial x}} \right] +
\frac{\partial }{\partial z}\left[ {\lambda (T,x,z)\frac{\partial
T}{\partial z}} \right]} \right\},}
\\ \\
{\displaystyle\  z_0 < z < Z, \quad l < x < d }
\end{array} \eqno{(4)}
$$

Boundary conditions for mould walls represent the character of heat
exchange on each sight of wall:

$$
\begin{array}{c}
{\displaystyle\ \left. {\lambda (T,z)\frac{\partial T}{\partial x}}
\right|_{x = d} = \alpha _1 \left( {T_{water} (\tau ,z) - \left. T
\right|_{x = d} } \right), \quad z_0 \le z \le Z,}
\\ \\
{\displaystyle\ \left. {\lambda (T,x)\frac{\partial T}{\partial z}}
\right|_{z = Z} = \alpha _2 \left( {T_{os.2} - \left. T \right|_{z =
Z} } \right), \quad l \le x \le d, \quad z = Z, }
\\ \\
{\displaystyle\ - \left. {\lambda (T,x)\frac{\partial T}{\partial
z}} \right|_{z = z_0 } = \alpha _3 \left( {T_{os.3} - \left. T
\right|_{z = z_0 } } \right), \quad l \le x \le d, \quad z = z_0, }
\\ \\
{\displaystyle\ \lambda (T,z)\left. {\frac{\partial T}{\partial x}}
\right|_{x = l + \delta } = \qquad \qquad \qquad \qquad \qquad
\qquad \qquad \qquad \qquad \qquad \qquad}
\\ \\
{\displaystyle\ \qquad \qquad \qquad = \frac{\lambda _{gz} }{\delta
}\left( { \left. T\right|_{x = l + \delta } - \left. T\right|_{x =
l} } \right) + \sigma_n\left[ {\left( {\frac{\left. T\right|_{x = l
+ \delta } }{100}} \right)^4 - \left( {\frac{\left. T\right|_{x = l}
}{100}} \right)^4} \right],}
\\ \\
\quad 0 \le z \le Z, \quad x = l + \delta ,
\\ \\
{\displaystyle\ \left. { - \lambda (T,z)\frac{\partial T}{\partial
x}} \right|_{x = l + \delta } = \alpha _4 \left( {T_{os.1} - \left.
T\right|_{x = d} } \right) + C_n\left[ {\left( {\frac{T_{os.1}
}{100}} \right)^4 - \left( {\frac{\left. T\right|_{x = d} }{100}}
\right)^4} \right], }
\\ \\
\quad z_0 \le z \le 0, \quad x = l + \delta,
\end{array} \eqno{(5)}
$$

\noindent where $d$ -- mold wall thickness, $z_0 $ -- mold wall
altitude over meniscus level, $\alpha _1 $ -- heat transfer
coefficient from the mould wall to cooling water, $T_{water} (\tau
,z)$ -- cooling water temperature in the mold channel, $\alpha
_{2,3,4} $ -- heat transfer coefficients from other mould wall to
environment, $T_{os.2,3,4} $ -- environment temperature, $C_n$ --
the resulted radiation coefficient.

The following balance equation describes distribution of cooling
water temperature in the mold channel:

$$
 c \cdot S \cdot v_{water} \frac{\partial T_{water} (\tau
,z)}{\partial z} = P_I \alpha _1 \left( {T_{water} (\tau ,z) -
\left. T \right|_{x = d} } \right) - P_E \alpha _E \left( {T_{water}
(\tau ,z) - T_E } \right), \eqno(6)
$$

\noindent where $c$ -- volume heat capacity of water, $S$ -- the
cross-section area of the mold channel, $v_{water} $ -- water
velocity, $P_I $ -- perimeter of the interior mold wall, $P_E $ --
perimeter of the external mold wall, $\alpha _E $ -- heat transfer
coefficient from cooling water to the external mould wall, $T_E$ --
external mould wall temperature.

The cooling water temperature on the entry in the mould channel is
known:

$$
 T_{water} (0,Z) = T_{water1} (\tau ) \eqno(7)
$$

\noindent and it's initial distribution in the mold channel:

$$
 T_{water} (0,z) = T_{water0} (z) \eqno(8)
$$

The following equation describes heat and mass transfer on the
curvilinear sections of CCM:

$$
\begin{array}{l}
{\displaystyle\ \frac{\partial T}{\partial \tau } + \theta _m (\tau
)\frac{\partial T(\tau ,r,\varphi )}{\partial \varphi } =
\frac{1}{c(T,r,\varphi )\rho (T,r,\varphi )} \times}
\\ \\
{\displaystyle\ \times \left\{ {\frac{\partial }{\partial r}\left(
{\lambda (T,r,\varphi )\frac{\partial T}{\partial r}} \right) +
\frac{1}{r^2} \cdot \frac{\partial }{\partial \varphi }\left(
{\lambda (T,r,\varphi )\frac{\partial T}{\partial \varphi }} \right)
+ \frac{\lambda (T,r,\varphi )}{r} \cdot \frac{\partial T}{\partial
r}} \right\} }
\end{array} \eqno{(9)}
$$

\noindent where $\theta _m $ -- angular velocity of ingot driving on
the $m$-th curvilinear section.

The conditions for unknown boundary on the curvilinear sections are

$$
\begin{array}{c}
{\displaystyle\ \left. {T(\tau ,r,\varphi )} \right|_{r = \xi _{1,2
- } (\tau ,\varphi )} = \left. {T(\tau ,r,\varphi )} \right|_{r =
\xi _{1,2 + } (\tau ,\varphi )} = T_{kr} ,}
\\ \\
{\displaystyle\ \left. {\lambda (T,r,\varphi )\frac{\partial
T}{\partial \bar {n}}} \right|_{\xi _{1 - } } - \lambda (T,r,\varphi
)\left. {\frac{\partial T}{\partial \bar {n}}} \right|_{\xi _{1 + }
} = \mu \rho _{kr} \left( {\theta _m (\tau ) \cdot \frac{\partial
\xi _1 }{\partial \varphi } + \frac{\partial \xi _1 }{\partial \tau
}} \right),}
\\ \\
{\displaystyle\ \quad \xi _1 (0,\varphi ) = \xi _{1_0 } (\varphi ),}
\\ \\
{\displaystyle\ \lambda (T,r,\varphi )\left. {\frac{\partial
T}{\partial \bar {n}}} \right|_{\xi _{2 + } } - \lambda (T,r,\varphi
)\left. {\frac{\partial T}{\partial \bar {n}}} \right|_{\xi _{2 - }
} = - \mu \rho _{kr} \left( {\theta _m (\tau ) \cdot \frac{\partial
\xi _2 }{\partial \varphi } + \frac{\partial \xi _2 }{\partial \tau
}} \right),}
\\ \\
{\displaystyle\ \quad \xi _2 (0,\varphi ) = \xi _{2_0 } (\varphi ),}
\end{array} \eqno{(10)}
$$

\noindent where $\xi _1 (\varphi )$ and $\xi _2 (\varphi )$ -- phase
boundaries (interfaces).

The boundary conditions of the secondary cooling zone include
radiant and convective components of heat exchange and account for
the complex mechanism of heat-conducting due to air-mist cooling
using compressed air and water. The boundary conditions on the
curvilinear sections are

$$
\begin{array}{l}
{\displaystyle\ - \left. {\lambda (T,\varphi )\frac{\partial
T}{\partial r}} \right|_{r = r_m } = \alpha _I (G_m (\tau ),\varphi
) \cdot \left( {T_{I _m} \left. { - T} \right|_{r = r_m } } \right)
+ C_{I _m} \left( {T_{I _m}^4 - \left. {(T} \right|_{r = r_m } )^4}
\right)}
\end{array} \eqno{(11)}
$$

$$
\begin{array}{l}
{\displaystyle\ \left. {\lambda (T,\varphi )\frac{\partial T_2
}{\partial r}} \right|_{r = r_m + 2l} = }
\\ \\
{\displaystyle\ = \alpha _E (G_m (\tau ),\varphi ) \cdot \left(
{T_{E _m} - \left. T \right|_{r = r_m + 2l} } \right) + C_{E _m}
\left( {T_{E _m}^4 - \left. {(T} \right|_{r = r_m + 2l} )^4}
\right),}
\end{array} \eqno{(12)}
$$

\noindent where $\alpha _I (G_m (\tau ),\varphi )$, $ \alpha _E (G_m
(\tau ),\varphi )$ -- convective heat transfer coefficients, $C_{I
_m}, C_{E _m}$  -- the resulted radiation coefficients, $T_{I _m} ,
 T_{E _m}$ -- environment temperatures, $G_m (\tau)$ -- water
discharge on the $m$-th section.

The following equation describes the heat and mass transfer on
rectilinear sections of CCM (analogously (1)):

$$
\begin{array}{l}
{\displaystyle\ \frac{\partial T}{\partial \tau } + v(\tau
)\frac{\partial T\left( {\tau ,x,z} \right)}{\partial x} = }
\\ \\
{\displaystyle\ = \frac{1}{c(T,x,z)\rho (T,x,z)}\left\{
{\frac{\partial }{\partial x}\left[ {\lambda (T,x,z)\frac{\partial
T}{\partial x}} \right] + \frac{\partial }{\partial z}\left[
{\lambda (T,x,z)\frac{\partial T}{\partial z}} \right]} \right\} }
\end{array} \eqno{(13)}
$$

When the liquid phase passes the straightening point on the
rectilinear section of the secondary cooling zone, the conditions
for the unknown phase boundary are set:

$$
\begin{array}{c}
{\displaystyle\ \left. {T(\tau ,x,z)} \right|_{x = \xi _{1,2 - }
(x,z)} = \left. {T(\tau ,x,z)} \right|_{x = \xi _{1,2 + } (x,z)} =
T_{kr} ,}
\\ \\
{\displaystyle\ \left. {\lambda (T,x,z)\frac{\partial T}{\partial
\bar {n}}} \right|_{\xi _{1 - } } - \left. {\lambda
(T,x,z)\frac{\partial T}{\partial \bar {n}}} \right|_{\xi _{1 + } }
= \mu \rho _{kr} \left( {v(\tau ) \cdot \frac{\partial \xi _1
}{\partial x} + \frac{\partial \xi _1 }{\partial \tau }} \right),}
\\ \\
{\displaystyle\ \left. {\lambda (T,x,z)\frac{\partial T}{\partial
\bar {n}}} \right|_{\xi _{2 + } } - \lambda (T,x,z)\left.
{\frac{\partial T}{\partial \bar {n}}} \right|_{\xi _{2 - } } = -
\mu \rho _{kr} \left( {v(\tau ) \cdot \frac{\partial \xi _2
}{\partial x} + \frac{\partial \xi _2 }{\partial \tau }} \right).}
\end{array} \eqno{(14)}
$$

The boundary conditions for the rectilinear section:

$$
\begin{array}{l}
{\displaystyle\ - \lambda (T,x)\left. {\frac{\partial T}{\partial
z}} \right|_{z = z_p } = \alpha _I (G_m (\tau ),x) \cdot \left( {T_I
- \left. T \right|_{z = z_p } } \right) + C_{I_4 } \left( {T_I ^4 -
\left. {(T} \right|_{z = z_p } )^4} \right)}
\\ \\
{\displaystyle\ \left. {\lambda (T,x)\frac{\partial T}{\partial z}}
\right|_{z = z_p + 2l} = }
\\ \\
{\displaystyle\ \qquad \qquad  =\alpha _E (G_m (\tau ),x) \cdot (T_E
- \left. T \right|_{z = z_p + 2l} ) + C_{E_4} \left( {T_E ^4 -
(\left. T \right|_{z = z_p + 2l} )^4} \right). }
\end{array} \eqno{(15)}
$$

We assume, that the thermal stream of the end of the rectilinear
site is equal to zero:

$$
 \left. {\lambda (T,z)\frac{\partial T}{\partial x}} \right|_{x =
x_f } = 0. \eqno(16)
$$

The initial conditions for all temperature field (on the rectilinear
and curvilinear sections):

$$
\begin{array}{l}
{\displaystyle\ T(0,x,z) = T_0 (x,z) }
\\ \\
{\displaystyle\ T(0,r,\varphi ) = T_0 (r,\varphi ). }
\end{array} \eqno{(17)}
$$

It is required to define the convective heat transfer coefficients
$\alpha _I (G_m (\tau ),\varphi )$, and  $\alpha _E (G_m (\tau
),\varphi )$ using the available information about ingot
temperature.

This is a boundary inverse problem and it is ill-posed in classical
sense. Well-posedness in classical sense (or Hadamard
well-posedness) means performance of three conditions: an existence
of a solution, its uniqueness and stability (input data continuous
dependence). In our case the third condition is not satisfied. This
is easily to verify using for the solution this problem the method
of direct reversion [2]. Therefore other approaches are necessary to
solve this problem.

\section{CHTC identification by least squares method} Consider an
ingot in first cooling section of secondary cooling zone. We have
ingot surface temperature measurements in some points. So we have to
solve the Dirichlet problem for interior heat exchange. The
finite-difference method was used to approximate the solution of
this problem. The convective heat-transfer coefficient (CHTC) has
special distribution along the surface of the ingot. Parabolic
function with a sufficient degree of accuracy approximates
distribution of CHTC on the part of surface that is exposed to
water-air spraying from one nozzle. This parabola has maximal value
in the point that corresponds to nozzle coordinate. CHTC is
considered as constant on the parts of the surface not subjected to
the forced cooling (fig. 2).

\begin{figure}[ht]
\centering
\includegraphics[scale=0.22, bb= 0 0 1246 594, clip]{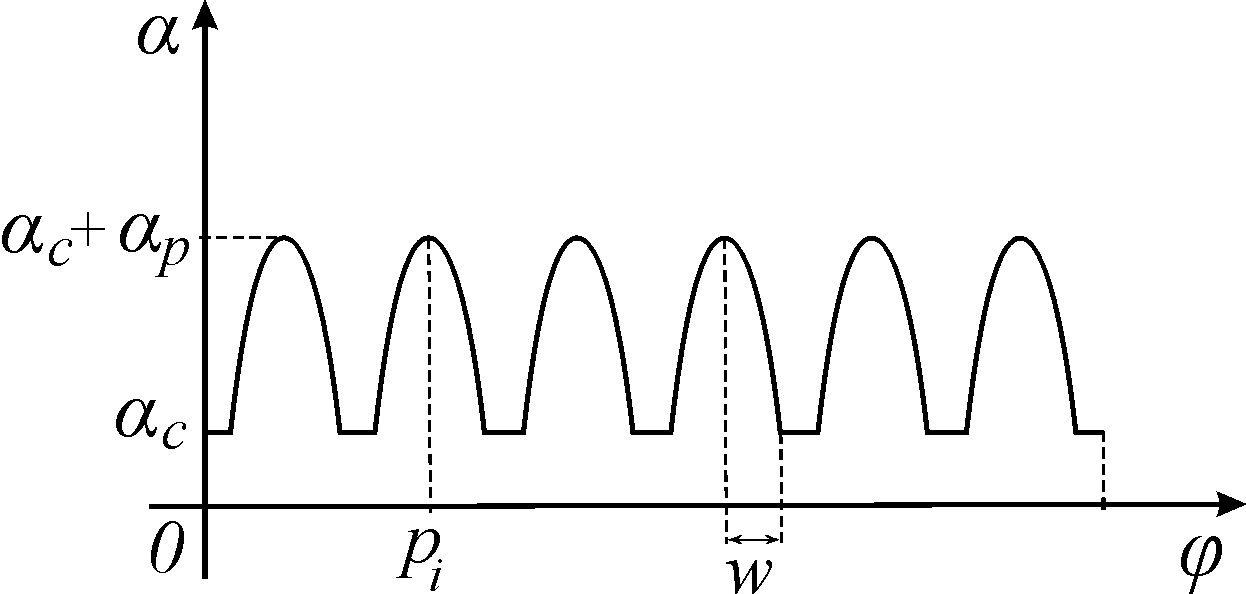}
\caption{}

\end{figure}

In one cooling section the same type spray nozzles are installed.
They give an identical water-air spray. Hence the CHTC is the same
parabola shifted along the abscissa axis (fig. 2).

All sites under spray nozzles can be reduced to the coordinate
origin so that the peak of each parabola should be over the
coordinate origin. Hence, it is necessary to define only two
parameters - $\alpha _p $ and $\alpha _c $. So, $\alpha (\varphi )$
is given by

$$
 \alpha (\varphi ) = \alpha _c - \frac{\alpha _p }{w^2}\varphi ^2 +
\alpha _p . \eqno (18)
$$

Consider the parts of the section, on which $\alpha (\varphi ) =
\alpha _c = const$. Let  $K$ be the ensemble of points $\varphi _i
$, in which CHTC is equal to constant. Let  $B$ be the ensemble of
other points.

The finite-difference approximation of boundary condition (11) is

$$
  \lambda _{i,0} \frac{T_{i,2} - 4T_{i,1} + 3T_{i,0}
}{2q} = \alpha _c (T_{I_1 } - T_{i,0} ) + C_{I_1 } (T_{I_1 } ^4 -
T_{i,0} ^4), \eqno (19)
$$

\noindent where $q$~-- step of finite-difference grid by radius $r_1
$ [3].

It follows that the discrepancy of heat flows on the boundary is:

\[
\Delta = \lambda _{i,0} \frac{T_{i,2} - 4T_{i,1} + 3T_{i,0} }{2q} -
C_{I_1 } \left( {T_{I_1 }^4 - T_{i,0} ^4} \right) - \alpha _c \left(
{T_{I_1 } - T_{i,0} } \right).
\]

Let us denote

\[
P_i = \lambda _{i,0} \frac{T_{i,2} - 4T_{i,1} + 3T_{i,0} }{2q} -
C_{I_1 } \left( {T_{I_1 }^4 - T_{i,0} ^4} \right), \quad Q_i =
T_{I_1 } - T_{i,0} .
\]

Then we find a value $\alpha _c $, such that the sum of squares of
discrepancies is minimum, i.e. the follow condition is satisfied

\[
S = \sum\limits_i {(P_i - \alpha _c Q_i )^2} \to \min , \quad \quad
\forall i:\varphi _i \in K.
\]

A necessary condition of the extremum existence of the function
$S(\alpha _c )$ is:

\[
\frac{\partial S}{\partial \alpha _c } = - 2\sum\limits_i {Q_i} (P_i
- \alpha _c Q_i ) = 0.
\]

It follows that

\[
\alpha _c = \frac{\sum\limits_i {Q_i P_i} } {\sum\limits_i {Q_i ^2}
}.
\]

To the each point $\varphi _i $ from В we will put in conformity a
point $y_i $ on the segment~$\left[ { - w,\,\,w} \right]$ such that
$\left| {y_i } \right|$ is equal to the distance from the
corresponding $\varphi _i $ to the coordinate of the nearest spray
nozzle. From (18) and (19) we gain a discrepancy

\[
\Delta = \lambda _{i,0} \frac{T_{i,2} - 4T_{i,1} + 3T_{i,0} }{2q} -
C_{I_1 } \left( {T_{I_1 }^4 - T_{i,0} ^4} \right) - \left( {\alpha
_c - \frac{\alpha _p }{w^2}y_i ^2 + \alpha _p } \right)\left(
{T_{I_1 } - T_{i,0} } \right).
\]

Then we can find a value $\alpha _p $, such that the sum

\[
S = \sum\limits_i {(P_i - (\alpha _c - \frac{\alpha _p }{w^2}y_i ^2
+ \alpha _p ) \cdot Q_i )^2} \to \min .
\]

From the following necessary condition  of extremum existence

\[
\frac{\partial S}{\partial \alpha _p } = 2\sum\limits_i {\left( {P_i
- \left( {\alpha _c - \alpha _p \left( {\frac{y_i ^2}{w^2} - 1}
\right)} \right)P_i } \right)\left( {Q_i \left( {\frac{y_i ^2}{w^2}
- 1} \right)} \right)} = 0
\]

\noindent we obtain $\alpha _p $

\[
{\displaystyle\ \alpha _p = \frac{ {\displaystyle\ \alpha _c
\sum\limits_i {Q_i ^2\left(  {\frac{y_i ^2}{w^2}}  - 1 \right)} }-
{\displaystyle\ \sum\limits_i {P_i Q_i \left( {\frac{y_i ^2}{w^2} -
1} \right)} } }{ {\displaystyle\ \sum\limits_i {Q_i ^2\left(
{\frac{y_i ^2}{w^2} - 1} \right)^2}} }. }
\]

On fig. 3 comparative results of calculations (1 -- by the method of
direct reversion, 2 -- by the least squares method) are presented.
For steel grade st40, width of a slab is 1m, l = 0,1m and v =
1(m/minute). The decision obtained by the method of direct reversion
is unstable and unsuitable for practical use. The second curve
represents a spline approximation, which is gained as a result of
the decision of a problem of identification by the least squares
method.

\begin{figure}[ht]
\centering
\includegraphics[scale=0.65, bb= 0 0 600 388, clip]{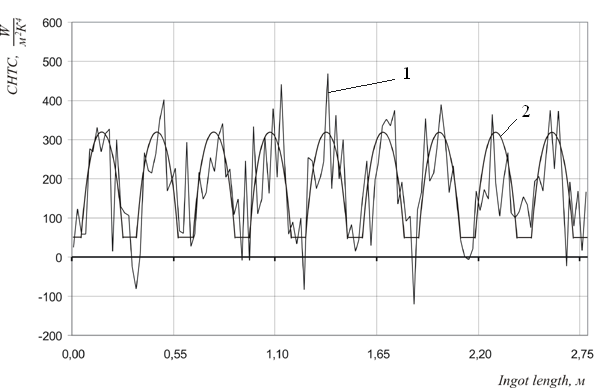}
\caption{}

\end{figure}

Thus, we fined the spline approximation of the CHTC, which is
distributed on the surface of the moving ingot. This approximation
gives the minimum of mean-square deviation between measured surface
temperature and calculated one according to the model as the result
of solving of the direct problem. The CHTC for other sections of the
secondary cooling zone is analogously defined. It should be noted
that an advantage of the offered method is that the estimation error
of the least squares method is negligibly small by relatively small
number of abnormal measurements. It is very important in case of
temperature measurement of a partially oxide scaled ingot surface.

\section {Operative adjustment of convective heat transfer
coefficient (CHTC)}

CHTC obtained by initial adjustment
 varies under changes of various parameters of process (for
example, ambient temperatures). Therefore, it is necessary to
provide its operative adaptation during work CCM. The fine-tuning of
parameters  should be carried out in real time. But during usual
work of CCM the information on a thermal condition of an ingot is
limited to temperature indications in small number of points of the
surface of an ingot. Such algorithms can be based on the stochastic
approximation method [4].


The temperature on the ingot surface is measured in every equal
small time intervals. Let us denote the measuring temperature data
$T_j^*$. The computer models the casting process using the presented
mathematical model. The under model calculated temperature in the
corresponding point we denote by $T_j$. It is necessary to correct
the model parameters using information about deviations between
measured and calculated temperature data to reduce these deviations
to minimum. The difficulty of the decision of the given problem is
that temperature measurements are deformed by a random telemetry
error.

Operative fine-tuning consists in refinement of the constant value
$\alpha _c $, which defines the distribution of the convective heat
transfer coefficient obtained by the solving of the problem of the
initial adjustment of parameters.


For using the algorithm of stochastic approximation it is necessary,
that the random error of temperature indications would have the zero
average and the finite variance.

The algorithm of parameter adjustment is

$$
 \alpha _{j + 1} = \alpha _j - k_j (T_j^\ast - T_j ), \eqno(20)
$$

\noindent where $\alpha _j $ -- $j$-th approximate value of $\alpha
_c $, $k_j $ -- special sequence of numbers, which satisfies to the
following conditions:

$$
\label{eq2} \mathop {\lim }\limits_{j \to \infty } k_j = 0, \quad
\sum\limits_{j = 1}^\infty {k_j = \infty } , \quad \sum\limits_{j =
1}^\infty {(k_j)^2 < \infty }. \eqno(21)
$$

For example the following elementary sequence satisfies to such
conditions

$$
k_j = \frac{a}{b + j},
$$

\noindent where $a,\;b\; \in R$, $a > 0$. Selecting numbers $a$ and
$b$, and also other sequences satisfying to the conditions (21), it
is possible to change speed of convergence of algorithm. In [3], for
example, it is recommended to keep $k_j$ as constant while the sign
of discrepancy $T_j^\ast - T_j $ not vary, and change then $k_j$ so
that to satisfy to above mentioned restrictions.

Truncation condition of the parameter fine tuning algorithm work is
occurrence of m last received approximations $\alpha_{n+1},
\alpha_{n+2},…,\alpha_{n+m}$ in a vicinity of $\alpha_n$ serves:

$$
|\alpha_n-\alpha_{n+i} |<\varepsilon, \quad   \forall i=1,…,m.
$$

If the condition is executed, assume $\alpha_c$ is equal $\alpha_n$.
For check we use values CHTC which have been picked up
experimentally at the decision of a direct problem of modeling of
thermal field CCM [1].

\section {Examples of realization of the stochastic approximation
method} Numerical modeling allows establishing the basic features of
trajectories of parameter fine-tuning process. On fig. 4
trajectories of parameter fine-tuning, characterizing a deviation of
the distributed parameter from true value, for the algorithm using
sequence $$k_j = \frac{a}{j}, \quad j = 1,2,3,...$$ are presented at
various values of factor $a$. When $a < 1$ very slow convergence is
observed. In this case the time of parameter tuning is inadmissible
big.

\begin{figure}[ht]
\centering
\includegraphics[scale=0.48, bb= 0 0 829 573, clip]{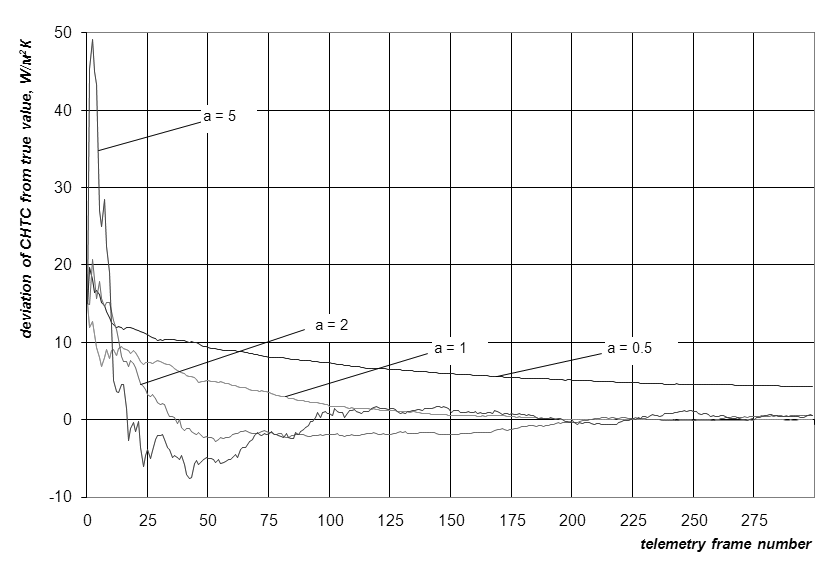}
\caption{}
\end{figure}

If to choose $a = 1$ the value of the parameter is in enough small
vicinity of true value approximately after 200th iteration. At $a =
2$ the trajectory of parameter fine-tuning reflects oscillations
with damped amplitude and frequency and not later than for 200
iterations the parameter is adjusted. At increase $a > 2$ the
amplitude of oscillations grows. In this case also oscillations with
damped amplitude and frequency are observed, but for fine-tuning it
is required considerably more iterations.

From here we conclude, that for the chosen sequence the best values
of the factor $a$ is a number from interval $1\leq a \leq 2$.

We investigate now influence of value $b$  on speed of the
algorithm's convergence. On fig.~5 trajectories of parameter
fine-tuning are shown for various values $b$. Values $b$ less than
zero lead to that fine-tuning go in a "wrong" direction while the
denominator is negative and at $i=-b$  the denominator is equal to
zero. Increase $b$  leads to decrease of a velocity of convergence
of algorithm. The same results have been obtained for sequences,
which will be described further. Therefore further parameter $b$
everywhere will be chosen equal to zero.

\begin{figure}[ht]
\centering
\includegraphics[scale=0.47, bb= 0 0 834 583, clip]{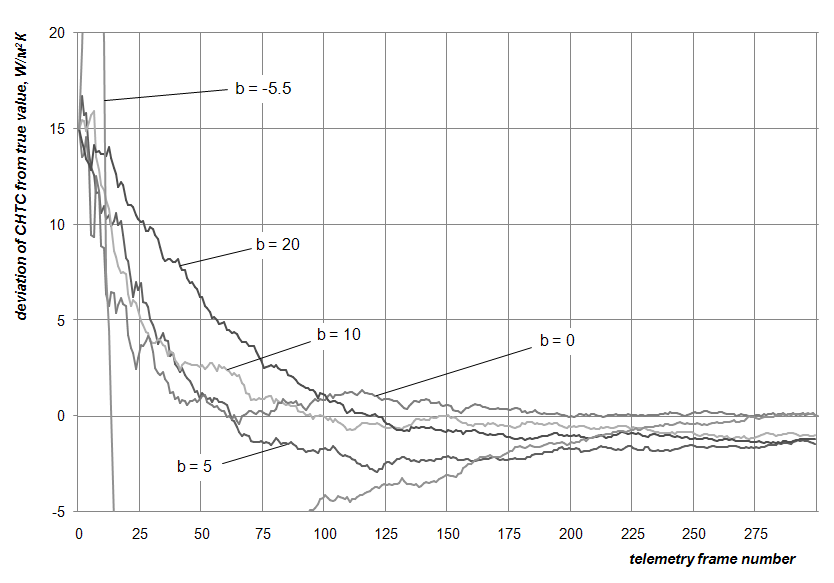}
\caption{}

\end{figure}

The following sequence also satisfies to conditions (21)

$$
 k_j = \frac{a}{n_j }, \quad n_{j + 1} = \left\{ {\begin{array}{l}
 n_j ,\quad \quad (T_j^\ast - T_j )(T_{j + 1}^\ast - T_{j + 1} ) > 0 \\
 j + 1,\;\quad (T_j^\ast - T_j )(T_{j + 1}^\ast - T_{j + 1} ) \le 0 \\
 \end{array}} \right.. \eqno (22)
$$

Results of this algorithm work are presented on fig. 6. In this case
factor $a$ needs to be chosen within $1 \le a \le 3$. Values out of
this range give smaller speed of algorithm convergence.

Consider another sequence, which also satisfies to conditions (21)

$$
 k_j = \frac{a}{n_j }, \quad n_{j + 1} = \left\{ {\begin{array}{l}
 n_j ,\quad \;\quad (T_j^\ast - T_j )(T_{j + 1}^\ast - T_{j + 1} ) > 0 \\
 n_j + 1,\quad \;(T_j^\ast - T_j )(T_{j + 1}^\ast - T_{j + 1} ) \le 0 \\
 \end{array}} \right.. \eqno (23)
$$

It has slower convergence than the previous two sequences. Results
of calculations with use of this sequence are presented on fig. 7.
Factor $a$ can be chosen within $0.5 \le a \le 2$. And, if $1.2 \le
a \le 1.5$, than obtained approximations differ from the true value
no more than on 6 {\%} after 20 iterations already.

\begin{figure}[ht]
\centering
\includegraphics[scale=0.6, bb= 0 0 650 505, clip]{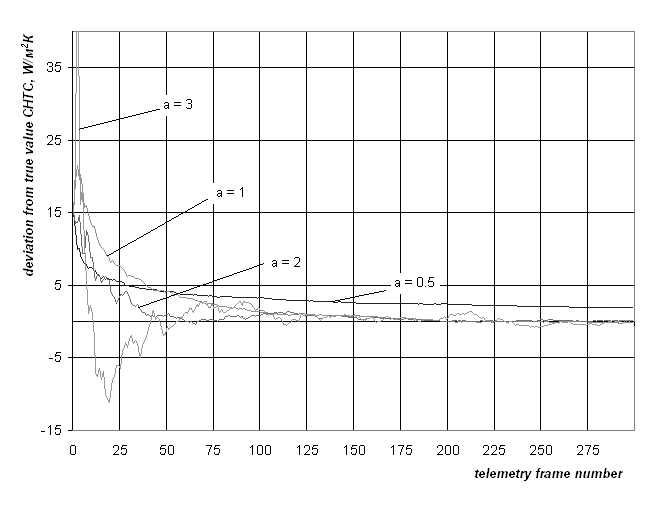}
\caption{}

\end{figure}

\begin{figure}[!ht]
\centering
\includegraphics[scale=0.5, bb= 0 0 770 585, clip]{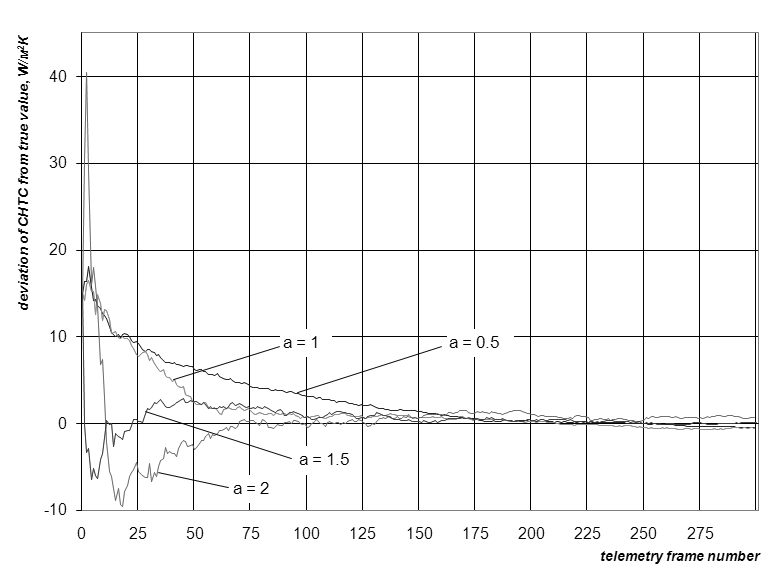}
\caption{}

\end{figure}
In the conclusion also it is necessary to notice, that the advantage
of stochastic approximation algorithm is its successful work for
enough wide interval of initial values of the distributed parameter.

\end{document}